\input amstex
\documentstyle{amsppt}
\magnification=1200
\pageheight{9 true in}

\def\o #1 {[#1]}
\def\so #1 {[{\ssize #1}]}
\def\wh #1{\widehat{#1}}
\def\wt #1{\widetilde{#1}}
\def\O{\Omega}

\def\tht{\thetag}

\def\k{\Bbbk}
\def\kb{\overline{\k}}

\def\e{\varepsilon}
\def\a{\alpha}
\def\b{\beta}
\def\g{\gamma}
\def\ep{\epsilon}
\def\Si{\Sigma}
\def\S{\Bbb S}
\def\A{\Bbb A}
\def\C{\Bbb C}
\def\Q{\Bbb Q}
\def\Z{\Bbb Z}
\def\Zp{\Bbb Z_{\ge 0}}
\def\ZZ{\Zp\times\Zp}
\def\nZ{\{0,\dots,n\}\times\Zp}
\def\itm#1{\item"(#1)"}
\def\l{\lambda}
\def\lto{\longrightarrow}
\def\x0n{x_0,\dots,x_n}
\def\1{\bar 1}
\def\ext{\operatorname{ext}}
\def\stab{\operatorname{stab}}
\def\we{\operatorname{weight}}
\def\const{\operatorname{const}}
\def\char{\operatorname{char}}
\def\Res{\operatorname{Res}}
\def\Eo{E${}_1$}
\def\Et{E${}_2$}
\def\lO #1{{}_{#1}\O}
\def\tP{\widetilde{P}}
\def\tO{\widetilde{\O}}
\def\LO{\Lambda^\O}

\topmatter
\title
On Newton interpolation of symmetric functions.
A characterization of interpolation Macdonald
polynomials.
\endtitle
\author
Andrei Okounkov
\endauthor
\address 
Department of Mathematics, University of
Chicago, 5734 South University Ave., Chicago, IL 60637-1546.
\endaddress
\email
okounkov\@math.uchicago.edu
\endemail
\thanks 
The  author is grateful to the
MSRI in Berkeley for hospitality and to the NSF for
financial support during his stay at the MSRI 
(grant DMS--9022140). 
\endthanks 
\toc\widestnumber\head{B.}
\head 1. Brief introduction \endhead
\head 2. General interpolation problem \endhead
\head 3. Examples of interpolation polynomials \endhead
\subhead 3.1 Universal interpolation polynomials \endsubhead 
\subhead 3.2 Factorial monomial symmetric functions \endsubhead 
\subhead 3.3 Factorial Schur functions \endsubhead 
\subhead 3.4 Interpolation Macdonald polynomials \endsubhead 
\head 4. Statement of the characterization theorem \endhead
\head 5. Reductions of the proof \endhead
\head 6. Proof of the theorem in the two variables case ($n=1$)\endhead
\head 7. Proof of the theorem for $n>1$ \endhead
\head A. Appendix. Table of perfect grids \endhead
\head B. Appendix. An 
example of the universal interpolation 
polynomial \endhead
\endtoc
\endtopmatter

\rightheadtext{Interpolation Macdonald polynomials}
\leftheadtext{Andrei Okounkov}

\document 

\head
1.~Brief introduction
\endhead

The purpose of this paper is to characterize
interpolation Macdonald polynomials inside a very general
Newton interpolation scheme for symmetric polynomials.
This general Newton interpolation problem is discussed
in Section 2; it depends as on a parameter on a map $\O$
$$
\left(\Zp\right)^{\text{\# of variables}}
@>\quad\O\quad>> \text{ ground field $\k$}\,,
$$
which we call a {\it grid} in $\k$. Our present understanding
of this general problem can be described as follows:
most of it is covered by an unexplored and mysterious 
ocean formed by generic grids $\O$ (see Section 3.1). In the
midst of this abyss there are 3 pieces of dry land, namely the
3 following exactly solvable cases:
\roster
\item factorial monomial symmetric functions,
\item factorial Schur functions,
\item interpolation Macdonald polynomials,
\endroster
which are described in Sections 3.2--3.4. The first case is
just dull, the second one is still rather elementary; both of
them are parameterized by an arbitrary sequence of pairwise distinct
elements of the ground field $\k$. These two continents 
are joined by a beautiful archipelago of interpolation 
Macdonald polynomials. More precisely, by interpolation 
Macdonald polynomials we mean the so called $BC$-type
interpolation Macdonald polynomials, introduced and 
studied in \cite{Ok4}. As particular cases and degenerations
these polynomials include polynomials studied by F.~Knop,
G.~Olshanski, S.~Sahi, and the author in a long series
of papers, see References. These polynomials depends on 5
parameters of which only 3 are non-trivial because of an
action of a 2-dimensional group of affine transformations. 
Some of the properties of these most remarkable polynomials are
discussed in Section 3.4.

It is natural to ask if any simple abstract property
characterizes the 3 above exactly solvable cases of our
general interpolation problem. As such a property we propose
the extra vanishing property \tht{4.2} which says that the
Newton interpolation polynomials should vanish not only at those
points where they are supposed to vanish by their definition
but also at certain extra points ``for free''. More precisely,
the polynomial labeled by a partition $\mu$ vanishes at the 
point labeled by a partition $\l$ unless $\mu$ is less or
equal to $\l$ in the partial order of partitions by inclusion
$\mu\subset\l$\,. We call all
grids that enjoy this property {\it perfect}. 

This extra vanishing property can be compared to the
following well-known property of ordinary Macdonald 
polynomials. Although the Gram-Schmidt orthogonalization
process requires a choice of a total order on the 
polynomials to be orthogonalized, the Macdonald
orthogonal polynomials do not actually depend on
the choice of a total order on the monomial symmetric
functions as long this total order is compatible with
the partial dominant order of partitions. 
\footnote{It would be probably interesting to
describe all interpolation or orthogonal 
symmetric polynomials which satisfy such a 
``extra triangularity'' condition.}

Also, the extra vanishing
property can be compared to a well known phenomenon in
integrable systems where many exactly solvable systems 
have ``extra'' integrals of motion, that is, more integrals
of motion than is required by the definition of integrability
\cite{Kr,CV}. It is interesting to notice that certain 
integrable many-body systems to which interpolation
Macdonald polynomials are very closely connected 
were conjectured to be characterized by this ``extra'' integrability
property \cite{CV}; later, however, certain new examples
were found in \cite{CFV}. 

Our situation is, of course, much simpler. Our main result
(Section 4) is that the three above cases plus degenerations
of the third one exhaust the set of all perfect grids. As
a corollary of this theorem, we conclude that no other grid
$\O$ admits a tableaux sum formula of the form \tht{4.1}
for the Newton interpolation
polynomials, nor does it admit an
integral representation of interpolation polynomials 
analogous to the $q$-integral representation for interpolation
Macdonald polynomials obtained in \cite{Ok4}. That is, any new
exactly solvable case of symmetric Newton interpolation has to
be based on some entirely new type of formulas.

The proof of this characterization theorem is given in
Sections 5--7. Section 5 contains some general statements,
whereas the two other sections are devoted to the
consideration of the many possible cases. 

There exist also non-symmetric Macdonald interpolation polynomials
which form a linear basis in the algebra of all polynomials,
see \cite{Kn,S2}. It is plausible that those polynomials might have
a similar characterization. Note also that a certain characterization
of ordinary Macdonald polynomials inside some general class of
orthogonal polynomials was found by S.~Kerov in \cite{K}.

\head
2.~General interpolation problem 
\endhead

We consider Newton interpolation of symmetric polynomials in 
$n+1$ variables ($n=0,1,\dots,\infty$) with coefficients in some
infinite field $\k$. First assume for simplicity that $n<\infty$;
the case $n=\infty$ will be covered at the end of this 
section. Any natural basis in the
space of symmetric polynomials of degree $\le d$ is indexed
by partitions
\footnote{ Our notation conventions about the partitions
are slightly different from the standard ones used in \cite{M1}.
For example, we number the parts of partitions as well as
coordinates of any vector starting from zero. This leads also
to a different definition of a diagram of a partition 
but makes certain formulas look more symmetric.}
$$
\l=(\l_0\ge\l_1\ge \dots \ge \l_n \ge 0)
$$
of length $\ell(\l)\le n+1$ (recall that $\ell(\l)$ is the number
of non-zero parts in $\l$) such that
$$
|\l|\le d \,.
$$
Here $|\l|=\l_0 + \dots +\l_n$. Therefore, the knots of our Newton
interpolation should be also indexed by partitions. In other words,
we need a function
$$
\mho : \big\{\l, \ell(\l)\le n+1\big\}  \lto \k^{n+1} \,,
$$
which takes a partition $\l$ to the corresponding knot of
interpolation $\mho(\l)$. We can construct such a function in the 
following way. Choose a function
$$
\O : \nZ \lto \k\,, \tag 2.1 
$$
which we shall call a {\it grid} in $\k$, and then set
$$
\mho(\l):=\left(\O(0,\l_0),\dots,\O(n,\l_n)\right) \in \k^{n+1} \,.
$$
To simplify notation, we shall use abbreviations
$$
\o ij :=\O(i,j)\,, \quad \wh \l  :=\mho(\l)\,.
$$
\definition{Definition 2.1} A grid $\O$ is said to be {\it non-degenerate}
if for any partition $\mu$ of length $\le n+1$ there exists a
symmetric polynomial
$$
P_\mu(x_0,\dots,x_n;\O) \in \k[\x0n]^{S(n+1)}
$$
satisfying the following Newton interpolation conditions
\roster
\item the degree of $P_\mu(x;\O)$ is $\le |\mu|$;
\item $P_\mu(\wh \l;\O)=0$ for all partitions $\l$ of length $\le n+1$
such that $|\l|\le|\mu|$ and $\l\ne\mu$;
\item $P_\mu(\wh \mu;\O)\ne0$.
\endroster
\enddefinition
\remark{Remark 2.2} It is clear that if the grid $\O$ is non-degenerate
then
\roster
\item all polynomials $P_\mu(x;\O)$ are uniquely defined up to
a non-zero factor from the field $\k$;
\item the polynomials $P_\mu(x;\O)$ as $\mu$ ranges over 
all partitions of length $\le n+1$ form a linear basis of the
vector space $\k[\x0n]^{S(n+1)}$ of all symmetric polynomials;
\item the degree of $P_\mu(x;\O)$ is precisely $|\mu|$.
\endroster
\endremark
The Newton interpolation polynomials $P_\mu(x;\O)$ generalize
polynomials considered by S.~Sahi in \cite{S1}. 

All non-degenerate grids admit a simple description.

\proclaim{Proposition 2.3} A grid $\O$ is non-degenerate if and
only if 
$$
\o ij \ne \o i'j' \,, \quad \forall i,j  \quad i\ge i', j<j' \,. \tag 2.2
$$
\endproclaim

Before we prove this proposition let us introduce the following
operation on grids \tht{2.1}.
For any $m=0,\dots,n$ we can define a grid
$$
\O_m: \{0,\dots,m\}\times\Zp \lto \k
$$
by simply restriction $\O_m(i,j):=\O(i,j)$ of the grid $\O$. 

Now the following proposition follows immediately from
the definition of polynomials $P_\mu(x;\O)$. (Since so far
these polynomials were defined up to a scalar factor only
the  equality \tht{2.3} is to be understood for the 
moment as an equality up to a scalar factor. However, below
in Definition 2.8 we shall choose a particular normalization 
of the polynomials
$P_\mu(x;\O)$ which will make \tht{2.3} into a
precise equality.)

\proclaim{Proposition 2.4} Suppose $\O$ is a non-degenerate grid
and $\mu$ is a partition such that $\mu_n=0$. Then
$$
P_\mu(x_0,\dots,x_{n-1},\o n0 ;\O)=P_{(\mu_0,\dots,\mu_{n-1})}
(x_0,\dots,x_{n-1};\O_{n-1}) \,. \tag 2.3 
$$
\endproclaim

\proclaim{Corollary 2.5} If a grid $\O$ is non-degenerate then so are
all grids $\O_m$, $m=0,\dots,n$.
\endproclaim

\demo{Proof of Proposition 2.3} Suppose that $\O$ is non-degenerate
and show that 
$$
\o ij \ne \o i'j' \,,\quad i\ge i', j<j' \,.
$$
By the above Corollary we can assume that $i=n$. Suppose that,
on the contrary, $\o nj = \o i'j' $ and $j'>j$. Then the 
following polynomial
$$
\prod_{m=0}^n \prod_{k=0}^{j'-1} (x_m - [n,k])
$$
is symmetric, has degree $n j'$ and vanishes at all points 
$\wh \l$ such that $|\l|\le nj'$. Thus, the ``only if'' part
of the proposition is established.

The ``if'' part of the proposition follows from the following
argument, which is an abstract form of the argument of S.~Sahi,
see \cite{S1}.

\proclaim{Proposition 2.6} Suppose a grid $\O$ satisfies the conditions
\tht{2.2} and let
$\S\subset\k$ be any ring such that
$$
[ij],\frac1{\o ij - \o i'j' }\in \S\,,\quad i\ge i', j<j' \,.
$$
Then for any $d=0,1,\dots$ and any function
$$
\phi : \big\{\l, |\l|\le d, \ell(\l)\le n+1\big\} \lto \S \tag 2.4
$$
there exists unique symmetric polynomial
$$
f_\phi \in \S[\x0n]^{S(n+1)}
$$
of degree $d$ such that
$$
f_\phi(\wh \l)=\phi(\l) \,. \tag  2.5
$$
\endproclaim
\demo{Proof} Induct on $n$ and $d$. The cases $n=0$ or $d=0$ are
clear. Suppose $n,d>0$. Introduce a linear map
$$
\ext : \S[x_0,\dots,x_{n-1}]^{S(n)} \lto \S[\x0n]^{S(n+1)}\,,
$$
which, by definition, takes a monomial symmetric function in
variables
$$
x_0-\o n0 ,\dots,x_{n-1}-\o n0 
$$
to the same monomial symmetric function in variables 
$$
x_0-\o n0 ,\dots,x_{n}-\o n0 \,.
$$
Observe that it is a degree preserving injection and that
$$
(\ext f)(x_0,\dots,x_{n-1},\o n0 )=f(x_0,\dots,x_{n-1})\,. \tag 2.6
$$
Now given $\phi$ we shall look for the solution  $f_\phi$ in
the form
$$
f_\phi=\ext f_1 + f_2 \, \prod_{i=0}^n (x_i - \o n0 )\,, \tag 2.7
$$
where $f_1$ and $f_2$ are unknown polynomials such that 
$$
\alignat2
f_1&\in \S[x_0,\dots,x_{n-1}]^{S(n)}\,,& \qquad \deg f_1 &= d\,, \\
f_2&\in \S[\x0n]^{S(n+1)}\,,& \qquad \deg f_2 &= d-n\,. 
\endalignat
$$
First, consider the equations \tht{2.5} for partitions $\l$ such that
$\l_n=0$. Since then the second summand in \tht{2.7} vanishes these
equations by \tht{2.6} and inductive assumption determine the polynomial
$f_1$. 

Now consider the equations \tht{2.5} for the remaining partitions $\l$
(that is, for $\l$ such that $\l_n>0$). Rewrite them in the form
$$
f_2(\wh \l)=\frac{\phi(\l) - (\ext f_1)(\wh \l)}
{\prod_i (\o i\l_i - \o n0 )}\,, \qquad \l_n>0\,, \tag 2.8
$$
and observe that the RHS of \tht{2.8} lies in $\S$. The set of
partitions $\l$ such that $|\l|\le d$ and $\l_n>0$ is in 
bijection
$$
\l \mapsto \l-\1=(\l_0-1,\dots,\l_n-1)
$$
with the set of partitions such that $|\l|\le d - n$. Therefore,
replacing $\O$ by $\O^1$ and using the inductive hypothesis we
can find a polynomial $f_2$ of degree $d-n$ satisfying the
equations \tht{2.8}. This proves existence of the required
polynomial $f_\phi$.

To verify uniqueness of $f_\phi$ it suffices to consider the
case $\S=\k$. We have an obvious $\k$-linear map from the
space of symmetric polynomials $f$ of degree $\le d$ to the
space of functions $\phi$ of the form \tht{2.4}, namely
$$
f \mapsto \phi \,,\quad  \phi(\l):=f(\wh \l)\,.
$$
Both spaces have same dimension and we just have proved that
this map is surjective. Hence it is an isomorphism with
the inverse map 
$$
\phi \mapsto f_\phi \,.
$$
This concludes the proof of the Propositions 2.6 and 2.3. \qed
\enddemo
\enddemo 

Let us introduce another operation on grids \tht{2.1}. 
Given a grid $\O$, we can for any 
$k=0,1,2,\dots$ define a new grid $\O^k$ 
$$
\O^k: \nZ \lto \k
$$
by the formula
$$
\O^k(i,j):=\O(i,j+k)\,.
$$
Then the following property follows immediately from
Definition 2.1 and Proposition 2.3 (As in the case of
Proposition 2.4, observe that 
since the normalization of the polynomials $P_\mu(x;\O)$
is yet to be specified 
the  equality \tht{2.9} is to be understood for the 
moment as an equality up to a scalar factor.) 

\proclaim{Proposition 2.7} Suppose $\O$ is a non-degenerate grid.
Then the grid $\O^1$ is also non-degenerate and for any
partition $\mu$ such that $\mu_n>0$ we have 
$$
P_\mu(x_0,\dots,x_n;\O)=P_{\mu-\1}(\x0n;\O^1) 
\, \prod_{i=0}^n (x_i - \o n0 )\,, \tag 2.9 
$$
where $\mu-\1$ stands for partition $(\mu_0-1,\dots,\mu_n-1)$. 
\endproclaim

Now suppose that $\O$ is a non-degenerate grid and
let us specify the normalization of the interpolation
polynomials $P_\mu(x;\O)$. Let us identify any partition $\l$
with its diagram which is, by definition, the following
subset of $\ZZ$
$$
\l=\big\{(i,j), j \le  \l_j -1 \big\} \subset \ZZ \,.
$$
Note that this definition differs by a 
coordinate shift by $1$ from the standard definition used in 
\cite{M1}. Let $\l'$ stand for the partition 
corresponding to the transposed diagram of $\l$, that is
$$
\l'_j:=\#\{ i, \l_i > j\}\,. \tag 2.10 
$$
We are free to normalize the polynomials $P_\mu(x;\O)$ by
setting their value at the point $\wh \mu$ to any non-zero
element of $\k$. We make the following choice

\definition{Definition 2.8} Suppose $\O$ is a non-degenerate grid and
$\mu$ is a partition with $\le n+1$ parts. Then we define
$P_\mu(x;\O)$ to be 
the unique polynomial satisfying the conditions of Definition 2.1
and the following normalization condition:
$$
P_\mu(\wh \mu;\O) = \prod_{(i,j)\in\mu} 
\big(\o i,\mu_i - \o \mu'_j-1,j \big)\,. \tag 2.11
$$
\enddefinition

The reason we prefer the normalization \tht{2.11} is the following 
\proclaim{Proposition 2.9} The normalization \tht{2.11} is the 
unique normalization compatible with \tht{2.3} and \tht{2.9}.
\endproclaim
\demo{Proof} Induct on $n$ and $|\mu|$. Since for
any partition $\mu$ either $\mu_n>0$ or $\mu_n=0$
we can apply either \tht{2.3} or \tht{2.9} to the
evaluation of 
$$
P_\mu(\wh \mu;\O)
$$
and thus reduce it to the case of smaller values of $|\mu|$
or $n$ respectively. \qed
\enddemo

By definition, let $\A\subset\k$ be the subring generated
by the elements
$$
\A=\Z\left[[ij],\frac1{\o ij - \o i'j' }\right]
\subset \k\,,\quad i\ge i', j<j' \,. \tag 2.12 
$$
Then since by \tht{2.11} we have
$$
P_\mu(\wh \mu;\O) \in \A
$$
for any $\mu$ we conclude from Proposition 2.6 and 
Definition 2.8 that 

\proclaim{Proposition 2.10}
$P_\mu(x;\O)\in\A[\x0n]^{S(n+1)}$\,.
\endproclaim  

\remark{Remark 2.11} Now let us explain the definition of
the polynomials $P_\mu(x;\O)$ in the case of infinitely
many variables ($n=\infty$). Suppose we are given a grid
$$
\O:\ZZ\lto \k \,. \tag 2.13
$$
Then  we need first a suitable
definition of the algebra of symmetric polynomials in
infinitely many variables. For any $m=0,1,2,\dots$
define a homomorphism
$$
\Res^{m+1}_m : \k[x_0,\dots,x_{m+1}]^{S(m+2)}
\lto  \k[x_0,\dots,x_{m}]^{S(m+1)}  \tag 2.14
$$
by the formula
$$
\left(\Res^{m+1}_m f\right)(x_0,\dots,x_{m})=
f(x_0,\dots,x_{m},\o m+1,0 )\,.
$$
Observe that
$$
\deg \Res^{m+1}_m f  \le \deg f \,.
$$ 
By definition, let $\LO$ be
$$
\LO:=\varprojlim \k[x_0,\dots,x_{m}]^{S(m+1)}
$$
the inverse limit of filtered (by the degree
of polynomials) algebras 
with respect to homomorphisms \tht{2.14}. Observe that
by construction for any
$$
f\in\LO
$$
its degree $\deg f$ is well defined and so are its
values at the points of the form
$$
f(\wh \l)=f\big(\o 0\l_0 ,\o 1\l_1 ,\o 2\l_2 ,\dots\big)\,,
$$
where $\l$ is a partition. 
Therefore the interpolation problem described in
Definition 2.1 is well defined for polynomials in
the algebra $\LO$, so it makes sense to ask
whether the grid $\O$ is non-degenerate and
what are the corresponding Newton interpolation
polynomials. Using definitions and Propositions
2.3 and 2.4 one easily checks that
\roster
\item a grid \tht{2.13} is non-degenerate if and only
if the conditions \tht{2.2} are satisfied (which is
also equivalent to non-degeneracy of all grids 
$\O_m$, $0\le m<\infty$);
\item for any partition $\mu$ the sequence
$$
\Big\{P_{(\mu_0,\dots,\mu_m)}(x_0,\dots,x_m;\O_m)
\Big\}_{m\ge\ell(\mu)}
$$
defines an element of $\LO$ which equals $P_\mu(x;\O)$.
\endroster 
\endremark

\head
3.~Examples of interpolation polynomials
\endhead

\subhead 3.1 Universal interpolation polynomials 
\endsubhead

Let $\k^u$ be the field
$$
\k^u=\Q(u_{ij})\,,\quad (i,j)\in\nZ\,,
$$
of rational functions in variables $u_{ij}$ and 
let the grid $\O^u$ be the grid
$$
\o ij ^u= u_{ij}\,,
$$
which is clearly non-degenerate (the superscript $u$ stands
here for ``universal''). We have the corresponding 
polynomials
$$
P_\mu(x;\O^u)\in\A^u[\x0n]^{S(n+1)}\,, \tag 3.1
$$
where the ring $\A^u$ is, according to the definition \tht{2.12},
the ring
$$
\A^u=\Z\left[u_{ij},\frac1{u_{ij} - u_{i'j'} }\right]
\subset \k^u\,,\quad i\ge i', j<j' \,.
$$
The polynomials \tht{3.1} are universal in the following
sense. For any field $\k$ and any non-degenerate grid
$$
\O : \nZ \lto \k
$$
we have a natural  homomorphism of rings
$$
\psi : \A^u \to \k\,,  \quad \psi(u_{ij})=\o ij  
$$
which is well defined by non-degeneracy of $\O$ 
and takes universal interpolation polynomials to 
the corresponding specific ones
$$
\psi \left(P_\mu(x;\O^u)\right)=P_\mu(x;\O)\,.
$$

One can compute a few universal polynomials explicitly
(see an Example in Appendix B) but 
the formulas get very soon very complicated. One is forced
therefore to look at least general but nicer examples.

The three other examples we shall consider in this section
will be {\it exactly solvable} in the following sense: there
exist rather simple closed formulas for the interpolation polynomials.
For simplicity we consider the case of finitely many variables
($n<\infty$). However, all  formulas are stable in the
sense described in Remark 2.11.

\subhead 3.2 Factorial monomial symmetric functions
\endsubhead

Suppose that
$$
\o ij = c_j\,,
$$
where $c_0,c_1,\dots$ are pairwise distinct elements of $\k$. Then it
is easy to see that the interpolation polynomials $P_\mu$ are
simply factorial monomial symmetric functions
$$
P_\mu(x;\O)=\frac1{\#\stab_{S(n+1)}\mu} \sum_{s\in S(n+1)}
\prod_{i=0}^n \,\,\prod_{j=0}^{\mu_i-1} \left(x_{s(i)} - c_j\right)\,,
\tag 3.2
$$
where $\#\stab_{S(n+1)}\mu$ is the number of permutations
$s\in S(n+1)$ which leave the vector $\mu=(\mu_0,\dots,\mu_n)$
invariant. It is also easy to see that the polynomials \tht{3.2} 
vanish at more more point then it is prescribed by their definition.
Namely, we have
$$
P_\mu(\wh \l;\O)=0\,, \quad\text{unless}\quad \mu\subset\l\,,
$$
where the notation $\mu\subset\l$ means that $\mu_i \le \l_i$
for all $i$ (which is equivalent to the diagram of $\mu$ being
a subset of the diagram of $\l$).

\subhead 3.3 Factorial Schur functions
\endsubhead

Suppose that
$$
\o ij = c_{j-i}\,,
$$
where $\dots,c_{-1},c_0,c_1,\dots$ are pairwise distinct elements of
$\k$. Introduce the following factorial Schur polynomial
(see \cite{M2} and references therein and also \cite{OO1})
$$
s_\mu(x;\O):= 
\frac{\det\big[ (x_i-c_{-n})\cdots(x_i-c_{\mu_j-j-1})
\big]_{0\le i,j\le n}}
{\prod_{0\le i < j\le n} (x_i - x_j)} \,. \tag 3.3
$$
First observe that \tht{3.3} is indeed a polynomial because the
numerator in \tht{3.3} is an anti-symmetric polynomial in $\x0n$
and hence is divisible by the denominator. It is also clear
that the ratio has degree $|\mu|$ and its top-degree component
is the classical Schur function $s_\mu(x)$.  

We shall check momentarily (borrowing the argument from \cite{Ok1},
Section 2.4; see also \cite{OO1})
that we have
$$
P_\mu(x;\O)=s_\mu(x;\O)\,. \tag 3.4
$$
First show that 
$$
s_\mu(\wh \l;\O)=0\,, \quad\text{unless}\quad \mu\subset\l\,. \tag 3.5
$$
Observe that the denominator in \tht{3.3} does not vanish at any point
of the form $\wh \l$. Therefore, it suffices to check that the
numerator vanishes unless $\mu\subset\l$. Suppose that $\l_k < \mu_k$
for some $k$. Then for all $0\le j\le k\le i \le n$ we have
$$
\l_i \le \l_k < \mu_k \le \mu_j 
$$
and the corresponding matrix element vanishes
$$
(c_{\l_i-i}-c_{-n})\cdots(c_{\l_i-i}-c_{\mu_j-j-1})=0 \,.
$$
Then the matrix 
$$
\big[ (c_{\l_i-i}-c_{-n})\cdots(c_{\l_i-i}-c_{\mu_j-j-1})
\big]_{0\le i,j\le n} \tag 3.6 
$$
has a block-triangular form with index sets $\{0,\dots,k-1\}$,
$\{k\}$, and $\{k+1,\dots,n\}$. The middle diagonal block is
zero and so is the determinant of \tht{3.6}. This proves \tht{3.5}.

Now let us compute $s_\mu(\wh \mu;\O)$. By the same argument
the matrix \tht{3.6} with $\l=\mu$ is triangular, therefore
$$
\align
s_\mu(\wh \mu;\O)&=\frac
{\prod_{i=0}^n (c_{\mu_i-i}-c_{-n}) \cdots
(c_{\mu_i-i}-c_{\mu_i-i-1})}
{\prod_{0\le i < j\le n} (c_{\mu_i-i} - c_{\mu_j-j})} \\
&=\prod_{(i,j)\in\mu} (c_{\mu_i-i} - c_{\mu'_j-j-1})\,,
\endalign
$$
where one goes from the first line to the second by a
standard combinatorial argument (used in e.g.\ proof of
the hook-length formula). This proves \tht{3.4}. 

There is a convenient tableaux sum formula (see \cite{M2} and 
\cite{GG}) for the polynomial \tht{3.3}. Given a partition $\mu$,
a {\it reverse} tableaux $T$ on $\mu$ is, by definition, a function
on the diagram of $\mu$
$$
\mu \owns (i,j) \overset T \to\longmapsto \{0,\dots,n\}  
\,,
$$
that decreases along the columns
$$
T(i,j) > T(i',j)\,, \quad  i > i'\,,
$$
and does not increase along the rows
$$
T(i,j) \ge  T(i,j')\,, \quad  j < j'\,.
$$  
With this notation the tableaux sum formula is
$$
s_\mu(x;\O)=\sum_T\, \prod_{(i,j)\in\mu} 
\left(x_{T(i,j)} - c_{j-i-T(i,j)}\right) \,, \tag 3.7
$$
where the summation ranges over all reverse tableaux $T$ 
on $\mu$. 

\subhead 3.4 Interpolation Macdonald polynomials
\endsubhead

These polynomials depend on 5 parameters $a,b,c,q,t$
from 
\footnote
{Actually, the parameters $a,b,c,q,t$ may lie in some 
extension of the field $\k$; then \tht{3.8} implies that
$q$ and $t$ satisfy quadratic equations with coefficients
in $\k$}
the field $\k$ such 
that $q,t\ne 0$ and the corresponding grid is given by
$$
\o ij = a + b q^j t^i + \frac c{q^j t^i} \in \k \,. \tag 3.8 
$$
The non-degeneracy of such a grid is assured by inequalities
$$
\alignat 2
q^k &\ne t^l\,,& \quad k&>0,\,l\ge 0\,,\\
b q^k t^l&\ne c\,,& \quad k&>0,\,l> 0\,,
\endalignat
$$
see more precise conditions in Appendix A. The simultaneous 
shift \& scale transformations
$$
x_i \mapsto C_1 x_i + C_2\,, \quad i=0,\dots,n,\, 
C_1\in\k\setminus 0,\,C_2\in\k\,,
$$
reduce the number of non-trivial parameters to following 3
$$
q,t,c/b\,,
$$
where the last one can also assume the value $\infty$. 

In full generality the corresponding interpolation polynomials
$P_\mu(x;\O)$ were considered in \cite{Ok4}.  Important particular cases and
degenerations of them were considered earlier by F.~Knop, G.~Olshanski,
S.~Sahi, and the author in a long series of papers, see 
References. 

The following tableaux sum formula for the polynomials $P_\mu(x;\O)$
follows by a change of variables from the statement of Theorem 5.2 in
\cite{Ok4}
$$
P_\mu(x;\O)=\sum_T\,\psi_T(q,t^{-1}) \prod_{(i,j)\in\mu} 
\left(x_{T(i,j)} - \o i+T(i,j),j \right) \,, \tag 3.9
$$
where the summation ranges over all reverse tableaux $T$ 
on $\mu$ and $\psi_T(q,t^{-1})$ is a certain two-parametric
weight of a tableaux $T$ introduced by I.~Macdonald in
\cite{M1}, section VI.7. This weight is a product of factors 
of the form
$$
(1-q^k t^l)\,, \quad  k,-l\in\Zp \,,
$$
and it appears in the tableaux sum formula for the 
ordinary Macdonald polynomials 
$$
P_\mu(x;q,t^{-1})=\sum_T\,\psi_T(q,t^{-1}) \prod_{(i,j)\in\mu} 
x_{T(i,j)}  \,. \tag 3.10 
$$
In \tht{3.10}, the LHS denotes the ordinary Macdonald polynomial
with parameters $q$ and $t^{-1}$ (we use Macdonald's notation
for it; it is not to be confused with our interpolation polynomials). 
Note that, in particular, \tht{3.10} is the highest degree term of 
\tht{3.9}. 

The 5-parametric grid \tht{3.8} admits many degenerations; they are 
all listed in Appendix A. As in the two previous examples, the
polynomials $P_\mu(x;\O)$ also enjoy the property
$$
P_\mu(\wh \l;\O)=0\,, \quad\text{unless}\quad \mu\subset\l\,,
$$
which can be seen from \tht{3.9} (see Lemma 4.1 below) but
is actually used in the proof of \tht{3.9}. 

The polynomials \tht{3.9} seem to be very distinguished
special polynomials. They have a wealth of applications
in various fields of mathematics, see e.g.\ \cite{Ok1,OO4-5,KOO,Ok5}. 

As to the practical interpolation, in the particular case when 
$$
bc=0
$$ there exists  a
nice efficient algorithm for Newton interpolation, see \cite{Ok3}.
It really speeds up for the
following (Jack) degeneration of the grid \tht{3.8}
$$
\o ij = \a + \b j +\b' i \,.
$$
A remarkable feature of that algorithm is that the coefficients
of the Newton interpolation expansion of any symmetric
polynomial can be found without actually computing these (rather
complicated) Newton interpolation polynomials.
That algorithm can be also used for expansion
in ordinary Macdonald and Jack polynomials. 

Note that in two particular cases
$$
t=1\quad\text{or}\quad t=1/q
$$
the number \tht{3.8} depends only on $j$ or 
$j-i$ respectively.  Hence for these values of $t$ 
the interpolation Macdonald polynomials 
become particular cases of 
factorial monomial symmetric functions
or factorial Schur functions. 

\head
4.~Statement of the characterization theorem
\endhead

The exactly solvable examples 2--4 of the previous section
have at least two following common features. The interpolation polynomials
admit a tableaux sum formula of the form
\footnote{
It is easy to see that the formula \tht{3.2} can be written
in the form \tht{4.1} with $\we(T)\in\{0,1\}$.
}
$$
P_\mu(x;\O)=\sum_T\,\we(T) \prod_{(i,j)\in\mu} 
\left(x_{T(i,j)} - \o i+T(i,j),j \right) \,, \tag 4.1
$$
where the summation ranges over all reverse tableaux $T$ 
on $\mu$ and 
$$
\we(T)\in\k
$$
is some weight of a tableaux $T$. The interpolation polynomials
also enjoy the following {\it extra vanishing} property
$$
P_\mu(\wh \l;\O)=0\,, \quad\text{unless}\quad \mu\subset\l\,. 
\tag 4.2
$$
In fact, \tht{4.2} is a consequence of \tht{4.1} as the following
argument (borrowed from \cite{Ok1}, Section 3.8) shows

\proclaim{Lemma 4.1} If $T$ is a reverse tableau on $\mu$ and
$\l$ is a partition then
$$
\prod_{(i,j)\in\mu} 
\left(\o T(i,j),\l_{T(i,j)} - \o i+T(i,j),j \right) = 0\,,
$$
unless $\mu\subset\l$.
\endproclaim

\demo{Proof} Suppose that 
$$
\o T(i,j),\l_{T(i,j)} - \o i+T(i,j),j \ne 0 \,, \quad \forall
(i,j)\in\mu\,.
$$
In particular, for $i=0$ we obtain
$$
\l_{T(0,0)}\ne 0\,,\l_{T(0,1)}\ne 1\,,\dots\,,
\l_{T(0,j)}\ne j\,, \dots \,. \tag 4.3 
$$
On the other hand, since $T$ is a reverse tableau we have
$$
\l_{T(0,0)} \le \l_{T(0,1)} \le \dots \,. \tag 4.4
$$
The inequalities \tht{4.3} and \tht{4.4} imply that
$$
\l_{T(0,j)} > j \,, \quad j=0,1,\dots \,.
$$
Again, since $T$ is a reverse tableaux we have
$$
T(0,j) < T(1,j) < \dots < T(\mu'_j-1,j) \tag  4.5 
$$
and also
$$
j < \l_{T(0,j)} \le 
\l_{T(1,j)}  \le  \dots \le \l_{T(\mu'_j-1,j)} \,. \tag  4.6 
$$
By the definition \tht{2.10} the inequalities \tht{4.5}
and \tht{4.6} yield that
$$
\l'_j \ge \mu'_j \,, \quad j=0,1,\dots \,,
$$
which is equivalent to $\mu\subset\l$\,. \qed
\enddemo

\remark{Remark 4.2} By employing the same argument as used in the proof 
of Theorem 5.1 in \cite{Ok4} one can show {\it a priori} that,
conversely, the extra vanishing \tht{4.2} implies existence of
a tableaux sum formula \tht{4.1}. We shall obtain this implication
as a corollary of our main theorem. 
\endremark 

\remark{Remark 4.3} It is also clear that the extra vanishing \tht{4.2}
follows immediately from any analog of the $q$-integral representation
established in \cite{Ok4} for the interpolation Macdonald polynomials in
the case $a,b,c,q,t\in\C$ and $|q|<1$.  
\endremark 

The above discussion justifies the following

\definition{Definition 4.4} We shall
call a  non-degenerate grid $\O$ 
{\it perfect} if the polynomials $P_\mu(x;\O)$ enjoy the
extra vanishing property \tht{4.2}.
\enddefinition

Our main result is the following 

\proclaim{Main Theorem} The following is the list of all perfect grids
$\O$:

\smallskip
\noindent
{\bf \Eo}. $\o ij = c_{j}$,
where $c_0,c_1,\dots$ are pairwise distinct elements of
$\k$; the corresponding interpolation polynomials are 
the factorial monomial symmetric function (see Section 3.2). 

\smallskip
\noindent
{\bf \Et}. $\o ij = c_{j-i}$,
where $\dots,c_{-1},c_0,c_1,\dots$ are pairwise distinct elements of
$\k$; the corresponding interpolation polynomials are 
the factorial Schur  function (see Section 3.3).

\smallskip
\noindent
{\bf I}. $\o ij = a+b q^j t^i + c q^{-j} t^{-i}$,
where $a,b,c,q,t$ are elements of a certain extension of
$\k$; the corresponding interpolation polynomials are 
the interpolation Macdonald polynomials (see Section 3.4).

\smallskip
\noindent
{\bf II-IV}. The grid $\O$ and the interpolation polynomials
are one of the degeneration of the previous case (see Appendix A).
\endproclaim

Since the two first cases are much simpler than the remaining ones we 
refer to them as to the 1st and 2nd elementary cases and number them
by \Eo\ and \Et. 

It follows from the above theorem together with Lemma 4.1 and
Remark 4.3 that if there exist any other
exactly solvable cases of the general interpolation problem
described in Section 2, then the formulas for the corresponding
interpolation polynomials should have some entirely new
structure. 

The proof of the theorem will be given in Sections 5--7.

\head
4.~Reductions of the proof 
\endhead

First from \tht{2.3} and \tht{2.9} one immediately derives the
following

\proclaim{Proposition 5.1} If a grid $\O$
$$
\O : \nZ \lto \k \tag 5.1
$$
is perfect then so are all grids $\O^k$, $k=0,1,\dots$, and
$\O_m$, $m=0,\dots,n$.
\endproclaim

Introduce one more operation on grids. Given a grid \tht{5.1} we can 
define a grid
$$
\lO l:\{0,\dots,n-l\} \lto \k\,, \quad l=0,\dots,n\,,
$$
by setting
$$
\lO l(i,j):=\O(i+l,j)\,.
$$
From Proposition 2.3 it is clear that the operation
$$
\O \longmapsto \lO l
$$
preserves non-degeneracy. We now plan to show that it preserves
perfectness  as well. First, we establish the following

\proclaim{Proposition 5.2} Let $\O$ be a perfect grid and let
$P_\mu(x;\O)$ be the corresponding Newton interpolation 
polynomials. We have
$$
P_\mu(x;\O)=x^\mu + \dots\,, \tag 5.2
$$
where $x^\mu=x_0^{\mu_0}\cdots x_n^{\mu_n}$ and 
dots stand for lower monomials in lexicographic order.
\endproclaim

\demo{Proof} First show that 
$$
P_\mu(x;\O)=c_\mu x^\mu + \dots\,, \quad c_\mu\in\k\,, \tag 5.3
$$
for certain constants $c_\mu$. Induct on $|\mu|$, the 
case $|\mu|=0$ being clear. Let $m$ 
$$
m:= \deg_{x_0} P_\mu(x;\O)
$$
be the degree of $P_\mu(x;\O)$ as of a polynomial in $x_0$.
If $m<\mu_0$ then  \tht{5.3} is established with $c_\mu=0$.

Suppose therefore that $m \ge |\mu_0|$ and consider the
leading coefficient 
$$
g(x_1,\dots,x_n):=\left[x_0^m\right] P_\mu(x;\O) \in \k[x_1,\dots,x_n]^{S(n)}  
$$
of $P_\mu(x;\O)$ as of a polynomial in $x_0$. By our
hypothesis $g$ 
is a non-zero polynomial of degree 
$$
\deg g \le |\mu|-\mu_0=\mu_1+\dots+\mu_n\,.
$$
We claim that 
$$
g\left(\o 1,\l_0 ,\dots, \o n\l_{n-1} \right)=0
$$
for all partitions $\l$ such that $\l_i < \mu_{i+1}$ 
for some $i=0,\dots,n-1$. Indeed, by the extra vanishing 
condition \tht{4.2} the polynomial
$$
P_\mu\left(x_0,\o 1,\l_0 ,\dots, \o n\l_{n-1} ;\O\right)
$$
has in this case infinitely many zeros
$$
x_0 = \o 0,\l_0 ,\o 0,\l_0+1 ,\o 0,\l_0+2 ,\dots 
$$
and hence vanishes identically. Since the grid $\lO 1$ in
non-degenerate we conclude that
$$
g=\const\, P_{(\mu_1,\dots,\mu_n)}(x_1,\dots,x_n;\lO 1)\,,
$$
which  by inductive assumption establishes \tht{5.3}.

Now recall that for any non-degenerate grid the polynomials
$P_\mu(x;\O)$ form a linear basis in the $\k$-linear space
of all symmetric polynomials. This immediately implies that
$$
\forall \mu \quad c_\mu \ne 0 \,.
$$
Let us renormalize the polynomials $P_\mu(x;\O)$.
Introduce new polynomials
$$
\tP_\mu(x;\O):=\frac 1{c_\mu} P_\mu(x;\O)
$$
for which we have 
$$
\tP_\mu(x;\O)=x^\mu + \dots\,.
$$
It is clear that these polynomials also satisfy the
equalities \tht{2.3} and \tht{2.9}. Therefore, by Proposition
2.9 we conclude that
$$
\tP_\mu(x;\O)=P_\mu(x;\O)\,,
$$
which concludes the proof of the proposition. \qed
\enddemo 

In fact, the above argument establishes more than just \tht{5.2}.
We have also proved the two following facts: 

\proclaim{Proposition 5.3} We have
$$
P_\mu(x;\O)=x_0^{\mu_0} P_{(\mu_1,\dots,\mu_n)}
(x_1,\dots,x_n;\lO 1) + \dots \,,
$$
where dots stand for terms of lower degree in $x_0$.
\endproclaim

\proclaim{Proposition 5.4} If a grid $\O$ is perfect then so 
are the grids $\lO 1, \lO 2,\dots$.
\endproclaim 

Let us combine the 3 operations 
$$
\O\mapsto\O_m\,,  \quad  
\O\mapsto\O^k\,,  \quad 
\O\mapsto\lO l
$$
as follows. Given a grid \tht{5.1} 
we define for all $0\le l \le m \le n$ and all
$k\ge 0$ a new grid 
$$
\lO l_m^k : \{0,\dots,m-l\}\times\Zp \lto \k
$$
by the formula
$$
\lO l_m^k(i,j):=\O(i+l,j+k)\,.
$$
Then by Propositions 5.1 and 5.4 we have

\proclaim{Proposition 5.5} If a grid $\O$
is perfect then so are all grids $\lO l^k_m$. 
\endproclaim

From now on we assume that $\O$ is a certain given perfect grid 
and our goal is to show that $\O$ is one of the grids listed in 
the Appendix A. 

Introduce the following notation. By the symbol $\e$ we
shall denote {\it some non-zero} element of $\k$
$$
\e \in \k\setminus 0 \,.
$$
The purpose of this notation is to denote irrelevant
overall factors in our equations. An example of an element
$\e$  is any product of factors of the form  
$$
\o ij - \o i'j' \,, \quad j<j',i\ge i'\,.
$$
Now we prove the following 

\proclaim{Proposition 5.6} Let $F$ be the following rational 
function 
$$
F(y_1,y_2,y_3,y_4,y_5):= 
\frac {y_4^2-y_4 y_1-y_4 y_5+y_2 y_3
+y_1 y_5-y_3^2}{y_2-y_3} \,. \tag 5.4
$$
Then we have 
$$
\o i+1,j+2 = 
F\big(\o i,j ,\o i+1,j ,\o i,j+1 ,\o i+1,j+1 ,\o i,j+2 
\big)\,, \tag 5.5
$$
for all $j$ and all $i\le n-1$. We also have 
$$
\o i+2,j+1 = 
F\big(\o i,j ,\o i,j+1 ,\o i+1,j , \o i+1,j+1 , \o i+2,j \big)\,, 
\tag 5.6
$$
for all $j$ and all $i\le n-2$. 
\endproclaim

\demo{Proof} Prove \tht{5.5}. 
By virtue of Proposition 5.5 it suffices to
establish \tht{5.5} for 
$$
n=1\,, \quad i=j=0 \,.
$$ 
Consider the condition
$$
P_{(3,0)}(\wh {(2,2)};\O)=0\,. \tag 5.7
$$
A direct computation shows that 
$$
\multline
P_{(3,0)}(\wh {(2,2)};\O)=\\
=\e(
\so 00 \so 02 -\so 00 \so 11 -{\so 01 }^{2}
+\so 01 \so 10 +\so 01 \so 12 -\so 02 \so 11 -
\so 10 \so 12 +{\so 11 }^{2}) \,,
\endmultline \tag 5.8
$$
where 
$$
\e=\frac{(\so 12 -\so 11 )(\so 12 - \so 10 )}{\so 01 - \so 10 }\,. 
$$
Since
$$
\so 01 - \so 10 \ne 0
$$
we deduce from \tht{5.7} that 
$$
\so 12 = F\big(\so 00 ,\so 10 ,\so 01 ,\so 11 ,\so 02 \big)\,. 
\tag 5.9 
$$
This proves \tht{5.5}.

Prove \tht{5.6}. By virtue of Proposition 5.5 it suffices to
establish it for 
$$
n=2\,, \quad i=j=0 \,.
$$ 
Consider the condition
$$
P_{(2,0,0)}(\wh {(1,1,1)};\O)=0\,.
$$
A direct computation shows that 
$$
\multline
P_{(2,0,0)}(\wh {(1,1,1)};\O)=\\
=\e(
\so 10 \so 21
-\so 01 \so 21 +\so 11 ^2-\so 11 \so 00 +\so 01 
\so 10 -\so 10 ^2 +\so 00 \so 20 -\so 11 \so 20 
)
\endmultline \tag 5.10 
$$
where 
$$
\e=\frac{\so 21 -\so 20 }{\so 10 - \so 01 }\,. 
$$
Therefore, 
$$
\so 21 = F\big(\so 00 ,\so 01 ,\so 10 ,\so 11 ,\so 20 \big)\,. 
\tag 5.11 
$$
This concludes the proof. \qed 
\enddemo

It is clear from the formula 
\tht{5.4} and non-degeneracy that the equalities
\tht{5.5} and \tht{5.6} can be reversed as follows

\proclaim{Corollary 5.7} We have 
$$
\o ij = F\big(\o i+1,j+2 , \o i,j+2 , \o i+1,j+1 ,
\o i,j+1 , \o i+1,j \big)\,, \tag 5.12
$$
for all $j$ and $i\le n-1$, and also 
$$
\o ij = F\big(\o i+2,j+1 , \o i+2,j , \o i+1,j+1 ,
\o i+1,j , \o i,j+1 \big)\,, \tag 5.13
$$
for all $j$ and $i\le n-2$.
\endproclaim

The above equalities immediately result in the following two
propositions. 

\proclaim{Proposition 5.8} Suppose that for some $i$ and $j$
we have
$$
\o i,j =\o i+1,j+1 \,. 
$$
Then $\O$ is a grid of type \Et\, that is, $\o i,j $
depends on $j-i$ only. 
\endproclaim
\demo{Proof}
Follows from \tht{5.5}, \tht{5.6}, \tht{5.12}, \tht{5.13},
and the following identity 
$$
F(z,u,v,z,w)=v\,. \qed  \tag 5.14
$$
\enddemo 

\proclaim{Proposition 5.9} Suppose that $n=1$ and for some $j$ 
we have 
$$
\o 0,j =\o 1,j \,, \quad \o 0,j+1 =\o 1,j+1 \,.
$$
Then $\O$ is a grid of type \Eo\, that is, for all $k\ge 0$
we have
$$
\o 0,k =\o 1,k \,. 
$$ 
\endproclaim
\demo{Proof}
Follows from \tht{5.5}, \tht{5.12}, and the following identity 
$$
F(u,u,v,v,w)=w \,. \qed \tag 5.15
$$
\enddemo

\head
6.~Proof of the Theorem in the two variables ($n=1$) case. 
\endhead

Consider the condition
$$
P_{(4,0)}(\wh {(3,2)};\O)=0
$$
and make the substitution \tht{5.11}. After that
substitution one computes:  
$$
P_{(4,0)}(\wh {(3,2)};\O)
=\e \, (\so 00 - \so 11 ) \,
F_1\big(\so 00 ,\so 10 ,\so 01 ,\so 11 ,\so 02 ,\so 03 \big)\,, 
\tag  6.1
$$
where $\e$ is the following invertible element
$$
\e = \frac{(\so 03 -\so 02 )(\so 12 -\so 11 )(\so 12 -\so 10 )}
{(\so 02 -\so 10 )(\so 10 - \so 01 )^2}
$$
and $F_1$ denotes the following polynomial 
$$
F_1(y_1,y_2,y_3,y_4,y_5,y_6):=y_6 G_1(y_2,y_3,y_4) -
G_0(y_1,y_2,y_3,y_4,y_5)\,,
$$
where
$$
\multline 
G_0(y_1,y_2,y_3,y_4,y_5):=
{y_1}{y_4}{y_2}-{y_1}{y_4}{y_5}-{y_1}{
y_2}{y_5}+{y_1}{y_5^2}+{y_3^2}{y_4}+
{y_3^2}{y_2}\\ -{y_3^2}{y_5}-{y_3}{y_4^2}-{y_3}{
y_4}{y_2}-{y_3}{y_2^2}+{y_3}{y_5^2}+
{y_4^2}{y_5}+{y_4}{y_2}{y_5}-{y_4}{y_5^2}+{
y_2^2}{y_5}-{y_2}{y_5^2}
\endmultline 
$$
and 
$$
G_1(y_2,y_3,y_4):= (y_3-y_2)(y_3-y_4)\,.
$$

Since the product \tht{6.1} vanishes 
we have the following alternatives. 

\medskip 
\noindent{\bf Case 1}. We have 
$\so 00 - \so 11 = 0$. 
Then by Corollary  we find ourselves in the \Et\ case.

\medskip 
\noindent{\bf Case 2}. We have 
$$
\so 03\, G_1\big(\so 10 ,\so 01 ,\so 11  \big) -
G_0\big(\so 00 ,\so 10 ,\so 01 ,\so 11 ,\so 02 \big)=0 \,. \tag 6.2
$$
Then, since the above equation 
is linear in $\so 03 $,  we have again two possibilities.  

\medskip 
\noindent{\bf Case 2.1}. We have 
$$
G_0\big(\so 00 ,\so 10 ,\so 01 ,\so 11 ,\so 02 \big) =
G_1\big(\so 10 ,\so 01 ,\so 11  \big) = 0 \,.  \tag 6.3
$$
Then the only solution of 
$$
G_1\big(\so 10 ,\so 01 ,\so 11  \big) = (\so 01 -\so 10 )(
\so 01 - \so 11 )=0
$$
compatible with the non-degeneracy is
$$
\so 01 = \so 11 \,. \tag 6.4  
$$
Substituting it into the first equation we obtain
$$
(\so 02 - \so 10 )(\so 02 - \so 01 )(\so 00 -\so 10 )=0 \,.
$$
From non-degeneracy we obtain
$$
\so 00 = \so 10 \,.
$$
Thus, by Proposition 5.9 we are in the \Eo\ case.  

\medskip 
\noindent{\bf Case 2.2}. We assume that the
 grid $\O$ is not of the
elementary types \Eo\ or \Et. Then we can uniquely determine 
$\so 03 $ from the values of 
$\so 00 $, $\so 10 $, $\so 01 $, $\so 11 $, and
$\so 02 $  by the formula  
$$
\so 03 = \frac{G_0\big(\so 00 ,\so 10 ,\so 01 ,\so 11 ,\so 02 \big)}
{G_1\big(\so 10 ,\so 01 ,\so 11  \big)} \,.
$$ 
Then using \tht{5.5}  we can also uniquely
determine $\so 13 $. We claim that we can 
determine $\so 04 $ as well. To that end we repeat the 
entire argument for the grid $\O^1$ which is also perfect
by Proposition 5.5.  We claim
that for the grid $\O^1$  the only possible case
is again the Case 2.2.  Indeed,  the grid $\O^1$ cannot be of types
\Eo\ or \Et\ because otherwise by Propositions 5.8, 5.9
 the grid $\O$ itself 
is of type \Eo\ or \Et. Therefore, we can 
uniquely determine $\so 04 $ and then
$\so 14 $ and then all the rest 
$$
\o ij \,, \quad i=0,1\,, \quad j\ge 5 \,.
$$
in the same manner. 

Now, this implies that our grid $\O$
is one of the types I--IV. Namely, we shall
check below in Proposition 6.1 that for any perfect grid $\O$
which is not of the
elementary types \Eo\ or \Et\ one can find
a grid $\tO$ of one of types I--IV such that
$$
\so 00 =\wt{\so 00 }\,, \quad 
\so 10 =\wt{\so 10 }\,, \quad 
\so 01 =\wt{\so 01 }\,, \quad 
\so 11 =\wt{\so 11 }\,, \quad 
\so 02 =\wt{\so 02 }\,. \quad 
$$
But then, since those five values uniquely determine
the rest, we conclude that $\O=\tO$.

Therefore, to finish the proof of the theorem in the
$n=1$ case it suffices to establish the following

\proclaim{Proposition 6.1} Denote by $\Si$ the following
subset of $\ZZ$
$$
\Si=\big\{(0,0),(1,0),(0,1),(1,1),(0,2),(0,3)\big\} \subset \ZZ \,.
$$
Then any perfect grid $\O$ belongs to one of
the following classes (compare with classification in Appendix A):

\smallskip\noindent
{\bf \Eo.} $\o 0j =\o 1j $, $j=0,1,2,\dots$.

\smallskip\noindent
{\bf \Et.} $\o 0,j =\o 1,j+1 $, $j=0,1,2,\dots$.

\smallskip\noindent
{\bf I.} There exist $a,b,c,q,t\in\kb$, $q,t\ne0,\pm1$ such that 
$$
\o ij = a+ b q^j t^i + \frac{c}{q^{j} t^{i}}\,, \quad \forall (i,j)\in\Si\,.
$$
The elements $a,b,c,q,t\in\kb$ are determined uniquely up to 
the following symmetry
$$
q\mapsto q^{-1}\,, \quad 
t\mapsto t^{-1}\,, \quad 
b\mapsto c\,, \quad 
c\mapsto b\,.
$$

\smallskip\noindent
{\bf II.} There exist unique $\a,\b,\b',\g\in\kb$ such that 
$$
\o ij = \a+ \b j +\b' i+\g(\b j +\b' i)^2 \,, \quad \forall (i,j)\in\Si\,.
$$

\smallskip\noindent
{\bf IIIa.} There exist unique $\a,\a',\b,\b'\in\k$ such that 
$$
\o ij = \a+ (-1)^j (\a'+\b j +\b' i)\,, \quad \forall (i,j)\in\Si\,.
$$

\smallskip\noindent
{\bf IIIb.} There exist $\a,\a',\b,\b'\in\k$ such that 
$$
\o ij = \a+ (-1)^i (\a'+\b j +\b' i)\,, \quad \forall (i,j)\in\Si\,.
$$
In this case only the numbers $\a+\a'$, $\a-\a'+\b'$, and $\b$ are
determined by the numbers $\o ij $, $(i,j)\in\Si$. 

\smallskip\noindent
{\bf IIIc.} There exist unique $\a,\a',\b,\b'\in\k$ such that 
$$
\o ij = \a+ (-1)^{i+j} (\a'+\b j +\b' i)\,, \quad \forall (i,j)\in\Si\,.
$$

\smallskip\noindent
{\bf IV.} There exist unique $\a,\b,\b',q\in\kb$ such that 
$$
\align 
\o 0j &= \a + \b q^j\,, \quad j=0,\dots,3\,, \\
\o 1j &= \a + \b' q^{-j}\,, \quad j=0,1\,.
\endalign  
$$
\endproclaim 

The proof will be based on the following lemma which can
be established by direct inspection

\proclaim{Lemma 6.2}
For any 4-tuple
$$
(w_0,w_1,w_2,w_3)\in\k^4
$$
satisfying 
$$
w_1 \ne w_2
$$
we have the 3 following mutually exclusive possibilities:
\roster
\itm i There exist $a,b,c,q\in\kb$, $q\ne 0,\pm1$ such that 
$$
w_j = a+ b q^j + \frac{c}{q^{j}} \,, \quad j=0,\dots,3\,.
$$
The elements $a,b,c,q\in\kb$ are determined uniquely up to 
the following symmetry
$$
q\mapsto q^{-1}\,, \quad 
b\mapsto c\,, \quad 
c\mapsto b\,.
$$
In this case
$$
w_0-3w_1+3w_2-w_3\ne 0\quad\text{and}\quad 
w_0+w_1-w_2-w_3\ne 0 \,.
$$
\itm{ii} $w_0-3w_1+3w_2-w_3=0$ and 
there exist unique $\a,\b,\g\in\k$ such that 
$$
w_j = \a+ \b j +\g\b^2 j^2 \,, \quad j=0,\dots,3\,.
$$
\itm{iii} $w_0+w_1-w_2-w_3=0$ 
and there exist unique $\a,\a',\b\in\k$ such that 
$$
\o ij = \a+ (-1)^j (\a'+\b j)\,, \quad j=0,\dots,3\,.
$$
\endroster
The last case is a subcase of second one if $\char\k=2$\,.
\endproclaim 

\demo{Proof of Proposition 6.1} Assume that we are not in the 
elementary  cases \Eo\ or  \Et\  and apply Lemma 6.2 to
the following 4-tuple
$$
\so 00 ,\so 01 ,\so 02 ,\so 03 \in\k \,.
$$
Recall that by non-degeneracy we have $\so 01 \ne\so 02 $. Therefore
we have 3 mutually exclusive possibilities which we shall
call cases {\bf (i)}, {\bf (ii)}, and  {\bf (iii)} respectively.

\smallskip\noindent
{\bf (i)} In this case we have
$$
\o 0j = a+b q^j+\frac{c}{q^j}\,, \quad j=0,\dots,3\,. \tag 6.5
$$
Assume for a moment that $bc\ne 0$. Choose  $t_0\in\kb$ so
that 
$$
\so 10 = a+bt_0+\frac c{t_0}\,. \tag 6.6
$$
Then the other root of this equation equals ${\dsize \frac{c}{bt_0}}$.
Since we exclude the elementary cases the
following equation is satisfied and is not identically zero
$$
\so 03\,  G_1\big(\so 10 ,\so 01 ,\so 11  \big) -
G_0\big(\so 00 ,\so 10 ,\so 01 ,\so 11 ,\so 02 \big)=0 \,. \tag 6.7
$$
Substituting \tht{6.5} and \tht{6.6} into \tht{6.7} we obtain
$$
\left(q^3 b - c\right)
\left(a+b q t_0+\frac{c}{q t_0}- \so 11 \right)
\left(a+\frac{b t_0}q+\frac{c q}{t_0}- \so 11 \right) = 0 \,.
$$
Since $\so 01 \ne\so 02 $ we have $b q^3 \ne c$. Hence either
$$
\so 11 = a+b q t_0+\frac{c}{q t_0}\,,
$$
in which case we are done by setting $t=t_0$, or
$$
\so 11 = a+\frac{b t_0}q+\frac{c q}{t_0}\,,
$$
in which case we set $t={\dsize \frac{c}{bt_0}}$. 

Now consider the case $bc=0$. We can assume that $c=0$. Set
$$
t=\frac{\so 10 -a}b\,.
$$
Then the equation \tht{6.7} implies that either
$$
\so 11 = a+b q t\,,
$$
which means we are in case I, or 
$$
\so 11 = a+\frac{b t}q\,,
$$
which brings us in the case IV.

\smallskip\noindent
{\bf (ii)} In this case we have
$$
\o 0j = \a+ \b j+\g \b^2 j^2\,, \quad j=0,\dots,3\,. \tag 6.8
$$
By non-degeneracy we have $\b\ne 0$. Again, assume for a moment
that $\g\ne 0$. Choose  $\b'_0\in\kb$ so
that 
$$
\so 10 = \a+\b'_0+\g {\b'_0}^2\,. \tag 6.9
$$
Then the other root of this equation equals $-\b'_0-\frac1{\g}$.
Solving the equation \tht{6.7} for $\so 11 $ we obtain that
either
$$
\so 11 = \a + \b + \b'_0 + \g( \b + \b'_0)^2\,,
$$
in which case we set $\b'=\b'_0$, or
$$
\so 11 = \a - \b + \b'_0 + \g( - \b + \b'_0)^2\,,
$$
in which case we set $\b'=-\b'_0-\frac1{\g}$ to obtain
$$
\so 11 = \a + \b + \b' + \g( \b + \b')^2\,.
$$

Now if $\g=0$ we set $\b'=\so 10 -\a$ and the two
possibilities for $\so 11 $ bring us in the cases
II and IIIb respectively.

Similarly, the consideration of the 
{\bf (iii)} case leads to cases
IIIa and IIIc. This concludes the proof. \qed
\enddemo

\head
7.~Proof of the theorem  for $n>1$
\endhead 

Suppose $\O$ is not of type \Et; then by Proposition 5.8
we can assume
$$
\so 00 \ne \so 11 \,.  \tag 7.1 
$$

Consider the condition
$$
P_{(3,0,0)}(\wh {(2,1,1)};\O)=0\,. \tag 7.2 
$$
Using the substitution \tht{5.11} one computes 
$$
\multline
P_{(3,0,0)}(\wh {(2,1,1)};\O)=\\
=\e (\so 11 - \so 00 )
\left(
\so 20 \, 
G_2\big(\so 00 ,\so 10 ,\so 01 ,\so 11 ,\so 02 \big) -
G_3\big(\so 00 ,\so 10 ,\so 01 ,\so 11 ,\so 02 \big)
\right)
\endmultline \tag 7.3 
$$
where 
$$
\e=\frac
{(\so 02 -\so 01 )(\so 21 -\so 20 )(\so 11 -\so 20 )\phantom{{}^2}}
{(\so 20 -\so 01 )(\so 10 -\so 02 )(\so 10 -\so 01 )^2}
$$
and $G_2$ and $G_3$ are the following polynomials
$$
\multline
G_2(y_1,y_2,y_3,y_4,y_5):=y_4 y_5 - y_4 y_3 + y_2 y_1 - y_2 y_3 
- y_1 y_5 + y_3^2\,,\\
G_3(y_1,y_2,y_3,y_4,y_5):=y_4^3 - y_4^2 y_2 - y_4^2 y_1 - y_4^2 y_3 - y_4 y_2^2 +
y_4 y_2 y_1+ y_4 y_2 y_5+\\
2 y_4 y_2 y_3 + y_4 y_1 y_3 - y_4 y_3^2 + y_2^3 - y_2^2 y_5 - y_2^2 y_3 +
y_2 y_5 y_3-y_2 y_3^2 - y_1 y_3 y_5 + y_3^3\,.
\endmultline \tag 7.4 
$$

Recall that we assume that $\so 00 \ne \so 11 $. Therefore, from \tht{7.2}
we have a linear equation  in  $\so 20 $ and hence two
possible cases:

\medskip 
\noindent{\bf Case 1}. We have
$$
G_2 = G_3 = 0 \,.
$$
Using \tht{7.1} we can express $\so 02 $ from $G_2=0$
$$
\so 02 =
\frac{\so 01 ^2 - \so 10 \so 01 - \so 11 \so 01 + \so 10 \so 00 }
{\so 00 -\so 11 }\,. \tag 7.5 
$$
Substituting this expression into the equation $G_3=0$ we obtain
$$
\e (\so 11 - \so 01 ) (\so 11 - \so 00 -\so 01 +\so 10 )
(\so 11 - \so 00 + \so 01 - \so 10 ) = 0\,, \tag 7.6
$$
where $\e=\so 11 - \so 10 $. Observe that it is impossible
to have $\so 11 =\so 01 $ because then from \tht{7.5} we 
conclude $\so 02 =\so 10 $, which contradicts non-degeneracy.
Similarly, it is impossible to have 
$$
\so 11 = \so 00 + \so 01 - \so 10
$$
because then from \tht{7.5}  we conclude
$$
\so 02 = \so 00 
$$
again in contradiction with non-degeneracy. Thus, the only
possible solution of \tht{7.6} and \tht{7.5} is
$$
\align
\so 11 &= \so 00 - \so 01 + \so 10 \,,\\
\so 02 &=2\so 01 - \so 00 \,.
\endalign
$$
That is, the grid $\O_2$ is of the type IIIb. Recall that
the parameters $\a,\a',\b,\b'$ are not uniquely determined
by the grid $\O_2$  (see Proposition 6.1). We can find such values
of $\a,\a',\b,\b'\in\k$ that
$$
\o ij = \a+ (-1)^i (\a'+\b j +\b' i)\,, \quad i+j\le 2\,.
$$
for all $i+j\le 2$. Then by \tht{5.6} we conclude
that 
$$
\o ij = \a+ (-1)^i (\a'+\b j +\b' i)\,, \quad \forall i\le 2
\quad \forall j\ge 0 \,. \tag 7.7
$$
Note that now the parameters $\a,\a',\b,\b'$ 
are  uniquely determined. If $n=2$
this finishes the consideration of the Case 1. 

If $n>3$ then we have to look at the equation
$$
P_{(2,1,0,0)}(\wh {(1,1,1,1)} )=0\,. \tag 7.8
$$
By \tht{7.7} one concludes from the identity
$$
F(u,u+z,v,v-z,w)=w+z
$$
and \tht{5.6} that 
$$
\so 31 = \so 30 - \b \,. \tag 7.9
$$
From \tht{7.7} and \tht{7.9} one computes that \tht{7.8} is
equivalent to
$$
\big(\a - \a' - 3\b' - \so 30 \big)
\big(\so 21 - \so 30 \big) = 0
$$
By non-degeneracy  we conclude that
$$
\so 30 = \a - \a' - 3\b' 
$$
and then by \tht{5.6} and \tht{7.7} it follows that
$$
\o 3j = \a - \a' - 3\b' - \b j \,.
$$
Thus, the grid $\O$ is of the type IIIb.
Repeating the above  argument for the grids
$\lO 1,\lO 2, \dots$ 
we conclude by Proposition 5.5
that the entire grid $\O$ is of the  type IIIb.

\medskip
\noindent
{\bf Case 2.} Now we assume that $\O$ is not
of the IIIb type and then we can determine $\so 20 $ from \tht{7.3}
$$
\so 20 = \frac
{G_3\big(\so 00 ,\so 10 ,\so 01 ,\so 11 ,\so 02 \big)}
{G_2\big(\so 00 ,\so 10 ,\so 01 ,\so 11 ,\so 02 \big)}\,. \tag 7.10
$$
After that using \tht{5.6} we can also determine
$\so 21 ,\so 22 ,\dots$. It is clear that if the 
grid $\O_1$ is of the type \Eo,I,II,IIIa,IIIc then
the grid $\O_2$ is of the same type with the same
parameters. If the grid
$\O_1$ is of type IV then 
$$
G_2 =0 \,, \quad G_3 \ne 0 \,,
$$
that is, the equation \tht{7.2} has no solutions. 
For $n=2$ this concludes the proof. 

If $n>2$ then we have to repeat the entire argument for the
grids $\lO 1,\lO 2, \dots$. In order to be able to do so
we have to verify that once the grid $\O_2$ is not of type 
IIIb then neither is the grid $\lO 1_3$. But this follows
from the just established classification of the perfect
grids for $n=2$. Namely, it follows that if $\O$ is a perfect
grid then 
$$
\left(\text{$\lO 1_2$ is of type IIIb}\right)
\Longrightarrow
\left(\text{$\O_2$ is of type IIIb}\right) \,.
$$
This concludes the proof of the theorem.

\head Appendix A. Table of perfect grids. 
\endhead

\noindent
{\bf \Eo.} The first elementary case
$$
\o ij = \g_j\,,
$$
where $\g_0,\g_1,\dots$ are arbitrary pairwise distinct. The
corresponding interpolation polynomials are factorial 
monomial symmetric functions. 

\medskip  
\noindent
{\bf \Et.} The second elementary case
$$
\o ij = \g_{j-i}\,,
$$
where $\dots,\g_{-1},\g_0,\g_1,\dots$ are arbitrary pairwise distinct.
The corresponding interpolation polynomials are factorial 
Schur functions. 

\medskip  
\noindent
{\bf I.} The generic case of Macdonald interpolation polynomials 
$$
\o ij = a+ b q^j t^i + \frac{c}{q^{j} t^{i}}\,,
$$
where $q,t\ne 0$. The non-degeneracy conditions are
$$
\alignat 2
q^k &\ne t^l \,,& \qquad &\forall k,l \quad k>0, 0\le l \le n\,,\\
b q^k t^l&\ne c \,,& \qquad &\forall k,l \quad k>0, 0\le l \le 2n\,.
\endalignat
$$
If $t=1$ or $t=1/q$ we hit the two elementary cases above.  

\medskip  
\noindent
{\bf II.} The case 
$$
\o ij = \a+ \b j +\b' i+\g(\b j +\b' i)^2
$$
can be obtained from I by setting
$$
\gather
a=\a-2\frac{\g}{h^2}\,, \quad b=\frac1{2h}+\frac{\g}{h^2}\,, \quad 
c=-\frac1{2h}+\frac{\g}{h^2}\,, \\
q=1+\b h\,, \quad t=1+\b' h
\endgather
$$
and letting $h\to 0$. The non-degeneracy conditions are
$$
\alignat 2
k\b &\ne l\b' \,,& \qquad &\forall k,l \quad k>0, 0\le l \le n\,,\\
\g(k\b+l\b')&\ne -1 \,,& \qquad &\forall k,l \quad k>0, 0\le l \le 2n\,.
\endalignat
$$

The cases II and III(abc) are impossible if $\char\k>0$\,.

\medskip  
\noindent
{\bf III(abc).} The case 
$$
\o ij = \a+ \ep^j{\ep'}^i(\a'+\b j +\b' i)
$$
can be obtained from I by setting
$$
\gather
a=\a\,, \quad b=\frac{\a'}2+\frac1{2h}\,, \quad 
c=\frac{\a'}2-\frac1{2h}\,, \\
q=\ep(1+\b h)\,, \quad t=\ep'(1+\b' h)
\endgather
$$
and letting $h\to 0$. Here $\ep,\ep'=\pm 1$; namely, we have three
subcases: {\bf (a)} $(\ep,\ep')=(-1,1)$, 
{\bf (b)} $(\ep,\ep')=(1,-1)$, {\bf (c)} $(\ep,\ep')=(-1,-1)$.

\medskip  
\noindent
{\bf IV.} This case exists only for $n=2$ when we can have
$$
\o 0j = \a + \b q^j\,, 
\quad \o 1j = \a + \b' q^{-j}\,.
$$
It is obtained from I by setting
$$
a=\a\,, \quad b=\b\,, \quad 
c=\b' t
$$
and then letting $t\to 0$. The non-degeneracy in this case is
equivalent to $\b,\b'\ne 0$ and $q^k\ne 1,\b'/\b$ for all $k>0$.  

\head
Appendix B. An example of a universal interpolation 
polynomial.
\endhead

Consider the case $n=1$ and
$$
\mu=(3,0)\,.
$$
Then the universal polynomial $P^u_\mu(x_0,x_1)$ lies in 
$$
P^u_\mu(x_0,x_1)\in
\Z\left[u_{00},u_{01},u_{02},u_{10},u_{11},
\frac1{u_{02}-u_{10}},\frac1{u_{01}-u_{10}}\right][x_0,x_1]\,.
$$
In the basis of the monomial symmetric functions $m_\l$ it is given
by the formula 
$$
P^u_\mu(x_0,x_1)=m_{30}+ 
\frac
{c_{21}(u)\, m_{21} + c_{20}(u)\, m_{20} + c_{11}(u)\, m_{11} +
c_{10}(u)\, m_{10} + c_{00}(u)}
{(u_{02}-u_{10})(u_{01}-u_{10})}\,,
$$
where $c_{21},\dots,c_{00}$ are the following polynomials 
$$
\align
c_{21}=&
u_{02}u_{00}-u_{02}u_{10}+u_{02}u_{01}- 
u_{02}u_{11}-u_{00}u_{10}+u_{00}u_{01}-u_{00}
u_{11}+u_{10}^2
\\&
-u_{10}u_{01}
+u_{10}u_{11}-
u_{01}u_{11}+u_{11}^2\,,\\
c_{20}=&
u_{02}^2u_{10}-u_{02}^2u_{01}-u_{02}u_{00}
u_{01}+u_{02}u_{10}u_{01}+u_{02}u_{10}
u_{11}-u_{02}u_{01}^2+u_{00}u_{10}u_{11}
\\&
-u_{10}^3
-u_{10}^2u_{11}+u_{10}u_{01}^2
+u_{10}
u_{01}u_{11}-u_{10}u_{11}^2\,,\\
c_{11}=&
-(u_{02}+u_{10}+u_{01}+u_{11}) (u_{02}
u_{00}-u_{02}u_{10}+u_{02}u_{01}-u_{02}
u_{11}-u_{00}u_{10}
\\&
+u_{00}u_{01}
-u_{00}u_{11}+u_{10}^2-u_{10}u_{01}+u_{10}u_{11}-u_{01}
u_{11}+u_{11}^2)\,,\\
c_{10}=&
u_{02}^2u_{00}u_{01}-u_{02}^2u_{10}^2-u_{02}^2
u_{10}u_{11}+u_{02}^2u_{01}^2+u_{02}
u_{00}u_{01}^2+u_{02}u_{10}^3-u_{02}u_{10}^2
u_{01}
\\&
-u_{02}u_{10}u_{01}u_{11}-u_{00}
u_{10}^2u_{11}-u_{00}u_{10}u_{11}^2+u_{10}^3
u_{01}+u_{10}^3u_{11}-u_{10}^2u_{01}^2
\\&
+u_{10}^2 u_{11}^2
-u_{10}u_{01}^2u_{11}+u_{10}u_{11}^3
\,,\\
c_{00}=&
-u_{02}^2u_{00}u_{01}^2+u_{02}^2u_{10}^2
u_{01}+u_{02}^2u_{10}^2u_{11}-u_{02}^2u_{10}
u_{01}^2-u_{02}u_{10}^3u_{01}-u_{02}u_{10}^3
u_{11}
\\&
+u_{02}u_{10}^2u_{01}^2+u_{02}u_{10}^2
u_{01}u_{11}+u_{00}u_{10}^2u_{11}^2-u_{10}^3
u_{01}u_{11}+u_{10}^2u_{01}^2u_{11}-u_{10}^2
u_{11}^3\,.
\endalign 
$$

\enddocument

\end